\documentclass[aap,secthm,MSNbibl,dvips]{arximspdf}
\usepackage{mathbh}
\usepackage{graphicx}

\doi{10.1214/08-AAP556}
\volume{19}
\issue{2}
\pubyear{2009}
\firstpage{641}
\lastpage{660}

\makeatletter
\newcommand{\tlong}{t_{\mathrm{long}}}
\newcommand{\Cyl}{\operatorname{Cyl}}
\newcommand{\Cyldeux}{\widetilde{Cyl}}

\newcommand\mincut{\operatorname{Mincut}}
\newcommand\Capflow{\operatorname{Capflow}}

\newcommand{\Cut}{\operatorname{Cut}}
\newcommand{\Int}{\operatorname{Int}}
\newcommand{\dis}{\operatorname{dis}}

\newcommand{\Ber}{\operatorname{Ber}}
\newcommand{\Div}{\operatorname{Div}}
\newcommand{\Divc}{\operatorname{div}}
\newcommand{\mumin}{\mu_{\min}}
\newcommand{\mumax}{\mu_{\max}}

\newcommand{\Dun}{C_1}
\newcommand{\Ddeux}{C_2}

\newcommand{\Z}{\mathbb{Z}}
\newcommand{\Zp}{\mathbb{Z}_{*}^2}
\newcommand{\Zd}{\mathbb{Z}^d}
\newcommand{\Zdeux}{\mathbb{Z}^2}
\newcommand{\R}{\mathbb{R}}
\newcommand{\Rd}{\mathbb{R}^d}
\newcommand{\Rdeux}{\mathbb{R}^2}
\renewcommand{\P}{\mathbb{P}}
\newcommand{\E}{\mathbb{E}}
\newcommand{\Ld}{\mathbb{L}^d}
\newcommand{\Ldeux}{\mathbb{L}^2}
\newcommand{\Ldeuxstar}{\mathbb{L}_{*}^2}
\newcommand{\Ed}{\mathbb{E}^d}
\newcommand{\Edeux}{\mathbb{E}^2}
\newcommand{\Edeuxstar}{\mathbb{E}^2_*}
\newcommand{\Edo}{\vec{\mathbb{E}}^d}
\newcommand{\Edod}{\overrightarrow{\mathbb{E}}^d}
\newcommand{\Edodeux}{\overrightarrow{\mathbb{E}}^2}

\renewcommand{\Ld}{\Ldeux}
\renewcommand{\Zd}{\Zdeux}
\renewcommand{\Ed}{\Edeux}
\renewcommand{\Rd}{\Rdeux}
\renewcommand{\Edod}{\Edodeux}

\newcommand{\as}{\mbox{a.s.}}

\renewcommand{\epsilon}{\varepsilon}
\renewcommand{\phi}{\varphi}
\renewcommand{\liminf}{\underline{\lim}}

\newcommand{\resp}{resp. }

\newcommand{\1}{\mathbh{1}}

\newtheorem{coro}[thm]{Corollary}
\newtheorem{prop}[thm]{Proposition}
\newtheorem{lemma}[thm]{Lemma}
\newproclaim{notation}{Notation}
\newproclaim{conjecture}[thm]{Conjecture}

\makeatother

\begin{document}
\begin{frontmatter}

\title{Capacitive flows on a 2D random net}
\runtitle{Capacitive flows on a 2D random net}

\begin{aug}
\author[A]{\fnms{Olivier} \snm{Garet}\ead[label=e1]{Olivier.Garet@iecn.u-nancy.fr}\ead
[label=u1,url]{http://www.iecn.u-nancy.fr/\texttildelow garet/}\corref{}}
\runauthor{O. Garet}
\affiliation{Institut \'Elie Cartan Nancy}

\address[A]{Universit\'e Henri Poincar\'e Nancy 1\\
Campus Scientifique, BP 239 \\
54506 Vandoeuvre-l\`es-Nancy Cedex\\
France\\
\printead{e1}\\
\printead{u1}}
\end{aug}

\received{\smonth{8} \syear{2006}}
\revised{\smonth{7} \syear{2008}}

%
\begin{abstract}
This paper concerns maximal flows on $\mathbb{Z}^2$ traveling
from a convex set to infinity, the flows being restricted by a random
capacity. For every compact convex set $A$, we prove that the maximal flow
$\Phi(nA)$ between $nA$ and infinity is such that
$\Phi(nA)/n$ almost surely converges to the integral of a deterministic
function over the boundary of $A$. The limit can also be interpreted
as the optimum of a deterministic continuous max-flow problem.
We derive some properties of the infinite cluster in supercritical Bernoulli
percolation.
\end{abstract}



\end{frontmatter}

\section{Introduction}

The problem of finding the maximum flow in a capacitive network is
undoubtedly the most known problem in the theory of operational research.
We know since Ford and Fulkerson that the search of a maximum
capacitive flow and that of the minimal cutset in a graph are
two sides of the same coin.
In the applications, the problems can come under a form or another.
Thus, this duality allows to choose the formulation which is the most
adapted to the mathematical treatment.

In the last decade, the min-cut formulation has been shown to be a practical
and useful tool for image segmentation (see Xiaodong~\cite{MR2296111}
or Estrada and Jepson~\cite{1069212}, for instance).
It is not surprising since image segmentation is precisely running the
scissors along the line of cut.
Let us assume for instance that we have a picture of a person
and that we want to cut around the face in such a way that the
background is rather white along the
break: if $\eta_x$ represents the blackness of the point $x$, then
one can try to minimize the ``cost'' $\sum_{x\in C} \eta_x$, where $C$
is a curve, which separates the face of the person (beforehand identified)
from the rest of the photograph.

We give in the present article a probabilistic
treatment of this kind of cutset problem: the darkness of the points is
given here by a collection of identically distributed random variables, and
we want to know to what extent the cost of the minimal cutset is determined
by the geometry of the form to be encircled.
If we reformulate the problem using the max-flow min-cut duality,
we have random capacities on the bonds of $\Zdeux$ and we study
the maximum flow that can be carried from the boundary of a given set
to infinity.
To be more specific, we fix a compact convex subset $A\subset\Rdeux$ and
study the asymptotic behavior of the maximal flow $\Phi(nA)$ between
$nA$ and infinity, which is also the cost of a minimal cutset
separating $nA$ from infinity.
We will see that the maximal flow $\Phi(nA)$ between $nA$ and infinity is
such that $\Phi(nA)/n$ almost surely converges to the integral of a
deterministic function over the boundary of $A$.

\section{Notation and results}

\subsection*{Flows}
Formally, let $\Edod=\{(x,y)\in\Zd\times\Zd\dvtx \|x-y\|_1=1\}$
and $\Ed=\{\{x,y\}\in\Zd\times\Zd\dvtx \|x-y\|_1=1\}$, where $\|\cdot\|_1$
is the $\ell^1$-norm: $\|(a,b)\|_1=|a|+|b|$.
As usual, we denote by $\Ld=(\Zd,\Ed)$ the unoriented square lattice.

We say that a map $f\dvtx\Edod\to\R$ is a flow if $f(x,y)=-f(y,x)$ holds for
each edge $(x,y)\in\Edod$.

Let $(t_e)_{e\in\Ed}$ be a family of positive numbers.

We say that $f$ is a capacitive flow from $A$ to infinity if it satisfies
%
\begin{equation}
\cases{
|f(x,y)|  \le t_{\{x,y\}}, &\quad for each bond $(x,y)\in\Edod$,\cr
\Div f (x)  =  0, & \quad for $x\in\Zd\backslash A$,
}
\end{equation}
where $\Div j(x)=\sum_{y\in\Zd; \|x-y\|_1=1} j(x,y)$.

We denote by $\Capflow(A,\infty)$ the set of capacitive flows from $A$
to infinity.
The aim is to study the maximal flow from a convex set $A$ to infinity,
that is,

\begin{equation}
\max\Biggl\{ \sum_{x\in A\cap\Zd}\Div j(x); j\in\Capflow(A,\infty) \Biggr\},
\end{equation}
%

when the $(t_e)_{e\in\Ed}$ are given by some collection of independent
identically distributed random
variables.

\subsection*{Links with first passage percolation}

The efficiency of methods coming from first passage percolation in studying
the maximum flow through a randomly capacitated network was initially
pointed out by Grimmett and Kesten~\cite{grimmett-kesten}: precisely,
they gave the asymptotic behavior of the maximum flow through the
bottom of a rectangle to its top as an application of their advances in
first-passage percolation.

As already mentioned, the point is the use of the max-flow min-cut
theorem~\cite{MR0079251}.
In the current setting, we can prove that
%
\begin{equation}
\label{maxflowmincutinfini}
\max\Biggl\{ \sum_{x\in A}\Div j(x); j\in\Capflow(A,\infty) \Biggr\}=\mincut
(A,\infty)\qquad\as,
\end{equation}
where $\mincut(A,\infty)$ is the minimum of the quantity $\sum_{e\in C}
t_e$, where $C$ is taken among the subsets of $\Ed$ that separate $A$
from infinity, or
more precisely that are such that every infinite path in $\Ld$ starting
from $A$ meets $C$. Such a set is called a cutset (relative to $A$).
With this definition, we can write
\begin{equation}
\quad \mincut(A,\infty)=\min\Biggl\{\sum_{e\in C} t_e; C\subset\Ed\mbox{ and }C
\mbox{ is a cutset relative to }A \Biggr\}.
\end{equation}

The cutsets of $\Ldeux$ can be characterized as follows: Let $\Zp
=\Zdeux
+(1/2,1/2)$, $\Edeuxstar=\{\{a,b\}; a,b\in\Zp\mbox{ and }\|a-b\|
_1=1\}$
and $\Ldeuxstar=(\Zp,\Edeuxstar)$. It is easy to see that
$\Ldeuxstar$
is isomorphic to $\Ldeux$.

For each bond $e=\{a,b\}$ of $\Ldeux$ (\resp$\Ldeuxstar$), let us
denote by
$s(e)$ the only subset $\{i,j\}$ of $\Zp$ (\resp$\Zdeux$) such that
the quadrangle $\mathit{aibj}$ is
a square in $\Rdeux$. $s$ is clearly an involution, and it is not
difficult to
see that $s$ is a one-to-one correspondence between the cutsets in
$\Ldeux$
and the sets in $\Edeuxstar$ that contain a closed path surrounding $A$.
If $C$ is minimal for inclusion, then $s(C)$ is just a path surrounding
$A$, so the quantity $\sum_{e\in e} t_e$ can be interpreted as the
length of the path
in a first-passage percolation setting on $\Zp$.

This leads us to recall a basic result in first-passage percolation:

Assume that $m$ is a probability measure on $[0, +\infty)$, such that

\begin{equation}
\label{leshypotheses}
m(0)<1/2 \quad\mbox{and}\quad
\exists c>0, \qquad \int_{[0, \infty)}\exp(c x) \,d m(x)< +\infty.
\end{equation}

Let $\Omega=[0,+\infty)^{\Ed}$ and consider the probability measure
$\P
=m^{\otimes\Ed}$ on $\Omega$. For $e\in\Ed$, we define $t_e(\omega
)=\omega_e$, thus the variables
$(t_e)_{e\in\Ed}$ are independent identically distributed random
variables with common law $m$.

For each $\gamma\subset\Ed$, we define $l(\gamma)=\sum_{e\in
\gamma} t_e$.
We denote by $d(a,b)$ the length of the shortest path from $a$ to $b$,
that is,
\[
d(a,b)=\inf\{l(\gamma);\gamma\mbox{ contains a path from }a\mbox{ to }b\}.
\]
%
Then by the Cox--Durrett shape theorem~\cite{coxdurrett}, there exists
a norm $\mu$ on $\Rd$ such that

\begin{equation}
\label{elementaire}
\lim_{\|x\|_1\to+\infty} \frac{d(0,x)}{\mu(x)}=1 \qquad\as
\end{equation}

We can also define $l^*$ by $l^*(A)=l(s(A))$ and a (random) distance
$d^*$ on $\Zp\times\Zp$ by
\[
d^*(a,b)=\inf\{l^*(\gamma);\gamma\mbox{ contains a path from
}a\mbox{
to }b\}.
\]

Since $\Ldeuxstar$ is isomorphic to $\Ldeux$, it is easy to see that
$d^*(\cdot,\cdot)$ enjoys the same asymptotic properties as $d(\cdot,\cdot)$ does.


\subsection*{Main results}

We first recall some common notation:
$\mathcal{H}^{1}$ is the $1$-dimen\-sional normalized Hausdorff measure,
$\lambda^2$ is the $2$-dimensional Lebesgue measure, $\Divc$ is the
usual divergence operator, and $C^1_c(\Rdeux,E)$ is the set of
compactly supported $C^1$ vector functions from $\Rd$ to $E$. Let
$A\subset\Rdeux$ be a Caccioppoli set. We denote by $\partial A$ its
boundary and by $\partial^* A$ its reduced boundary, that is
constituted by the points $x\in\partial A$, where $\partial A$ admits an
unique outer normal, which is denoted by $\nu_A(x)$.

The main goal of the paper is the following theorem.

\begin{thm}\label{main}
We suppose that $m(0)<1/2$ and that
\[
\int_{[0, \infty)}\exp(c x) \,d m(x)<+\infty
\]
for some $c>0$.
Then for each bounded convex set $A\subset\Rdeux$ with 0 in the
interior, we have
\[
\lim_{n\to+\infty} \frac{\mincut(nA,\infty)}{n}=\int_{\partial
^* A}\mu
(\nu_A(x)) \,d\mathcal{H}^{1}(x).
\]

Equivalently,
\begin{eqnarray*}
& &\lim_{n\to+\infty} \frac1{n}\max\Biggl\{ \sum_{x\in nA\cap\Zdeux
}\Div
j(x); j\in\Capflow(nA,\infty) \Biggr\}\\
&&\qquad =  \sup\biggl\{\int_{A} \Divc f
\,d\lambda^2(x);f\in C^1_c(\Rdeux,\mathcal{W}_{\mu}) \biggr\},
\end{eqnarray*}
where
\[
\mathcal{W}_{\mu}=\{x\in\Rdeux\dvtx \langle x,w\rangle\le\mu(w)\mbox{ for
all }w\}.
\]
\end{thm}

Note that $\mathcal{W}_{\mu}$ is sometimes called the Wulff crystal
associated to $\mu$.

If we observe the last equality, we can see that
the optimal value of a discrete random max-flow
problem converges (after a suitable renormalization) to
the optimum of a deterministic continuous max-flow problem.

In fact, we even have exponential bounds for the fluctuations around
%
\begin{equation}\label{def-I-A}
\mathcal{I}(A)=\sup\biggl\{\int_{A} \Divc f  \,d\lambda^2(x);f\in
C^1_c(\Rdeux
,\mathcal{W}_{\mu}) \biggr\}.
\end{equation}

Indeed, we prove the following theorem.
\begin{thm}
\label{vitesseexpo}
Under the assumptions of Theorem~\ref{main}, it holds that
for each $\epsilon>0$,
there exist constants $\Dun,\Ddeux>0$, depending on $\epsilon$ and $m$,
such that
%
\begin{equation}
\label{les-deux-cotes}
\forall n\ge0\qquad\P\biggl(\frac{\mincut(nA,\infty)}{n\mathcal
{I}(A)}\notin(1-\epsilon,1+\epsilon) \biggr)\le\Dun\exp(-\Ddeux n).
\end{equation}
\end{thm}

With the help of Menger's theorem, we obtain the following corollaries.
\begin{coro}\label{avecmenger}
We consider supercritical Bernoulli percolation on the square lattice,
where the edges are open with probability $p>p_c(2)=1/2$.
Then for each bounded convex set $A\subset\Rdeux$ with 0 in the
interior, the maximal number
$\dis(A)$ for a collection of disjoint open paths from $A$ to infinity
satisfies
%
\begin{equation}
\qquad\exists\Dun,\Ddeux>0\ \forall n\ge0\qquad\P\biggl(\frac{\dis
(nA)}{n\mathcal{I}(A)}\notin(1-\epsilon,1+\epsilon) \biggr)\le\Dun\exp
(-\Ddeux n),
\end{equation}
where $\mathcal{I}(A)$ is the quantity defined in (\ref{def-I-A}), the
law $m$ of passage times being the Bernoulli distribution $(1-p)\delta
_0+p\delta_1$.
\end{coro}

This corollary has itself an easy and pleasant consequence.

\begin{coro}
\label{avecmenger2}
We consider supercritical Bernoulli percolation on the square lattice.
For each integer $k$, there almost surely exist $k$ disjoint open
biinfinite paths.
\end{coro}

Note, however, that this amusing corollary is not really new; indeed,
it can be
obtained as a consequence of Grimmett and Marstrand~\cite{Grimmett-Marstrand}---see also Grimmett~\cite{grimmett-book}, page
148, Theorem 7.2.(a).

The paper is organized as follows. In Section~\ref{lapremiere}, we recall
some basic properties in first-passage percolation and prove some
useful properties of the functional $\mathcal{I}$. Next, the proof of
Theorem~\ref{vitesseexpo} naturally falls into two parts: Section~\ref
{ladeuxieme} deals with
the upper large deviations appearing in the Theorem, whereas
Section~\ref{latroisieme} is about the lower ones.
We complete the proof of Theorem~\ref{main} and establish the
corollaries in Section~\ref{laquatrieme}.
In the final section, we discuss the possibility of an extension to
higher dimensions.

\section{Preliminary results}
\label{lapremiere}

\begin{notation*}
We denote by $\langle\cdot,\cdot\rangle$ the natural scalar product
on $\Rd$ and by $\|\cdot\|_2$ the associated norm. $\mathcal S$ is the
Euclidean unit sphere: $\mathcal S=\{x\in\Rd\dvtx\|x\|_2=1\}$.
\end{notation*}

\subsection{First-passage percolation}

Let us introduce some notation and results related to first passage
percolation. As previously, we suppose that~(\ref{leshypotheses}) is
satisfied and write $\mu$ for the norm given by~(\ref{elementaire}).

It will be useful to use
\[
\mu_{\max}=\sup\{\mu(x);\|x\|_1=1\}\quad \mbox{and}\quad\mu_{\min}=\inf\{\mu(x);\|x\|_1=1\}.
\]
Of course, $0<\mumin\le\mumax<+\infty$ and we have
\[
\forall x\in\Rd\qquad\mumin\|x\|_1\le\mu(x)\le\mumax\|x\|_1.
\]

The speed of convergence in equation~(\ref{elementaire}) can be specified:

\begin{prop}[(Large deviations, Grimmett--Kesten~\cite{grimmett-kesten})]\label{mugs}
$\!$For each \mbox{$\epsilon>0$}, there exist $C_3,C_4>0$ such that
\[
\forall x\in\Zd\qquad\P\bigl( d(0,x)\in[(1-\epsilon)\mu(x),(1+\epsilon
)\mu(x)] \bigr) \ge1- C_3 \exp(-C_4 \|x\|_1).
\]
\end{prop}

Note that in \cite{grimmett-kesten}, the proof of this result is only written
in the direction of the first axis, that is, for $x=ne_1$.
Nevertheless, it applies in any direction and computations can be
followed in order to preserve a
uniform control, whatever direction one
considers. See, for instance, Garet and Marchand~\cite{GM-large} for a
detailed proof in an analogous situation. The control of $\P(
d(0,x)>(1+\epsilon)\mu(x) )$ could also be obtained as a byproduct of
the foregoing Lemma~\ref{GDlong}.

\subsection{Properties of $\mathcal{I}$}

Since $\mu$ is a norm, it is obviously a convex function that does not
vanish on the Euclidean sphere $\mathcal{S}$.
So, it follows from Proposition 14.3 in Cerf~\cite{cerf-stflour} that
the identity
%
\begin{equation}
\label{cerf}
\int_{\partial^* A}\mu(\nu_A(x))\, d\mathcal{H}^{1}(x)=\mathcal{I}(A)
\end{equation}
holds for every Cacciopoli set, and particularly for compact convex
sets and polygons.

From equation~(\ref{def-I-A}), it is easy to see that
%
\begin{equation}
\label{homogene}
\mathcal{I}(\lambda A)=\lambda\mathcal{I}(A)
\end{equation}
holds for each Borel set $A$ and each $\lambda>0$.

\begin{lemma}
$\mathcal{I}(A)>0$ for each convex set $A$ with nonempty interior.
\end{lemma}
\begin{pf}
For each $x\in\partial^* A$, $\mu(\nu_A(x))\ge\mumin\|\nu_A(x)\|
_1\ge
\mumin\|\nu_A(x)\|_2=\mumin$, so it follows from~(\ref{cerf}) that
$\mathcal{I}(A)\ge\mumin\mathcal{H}^1(\partial A)$.
\end{pf}

The next lemma clarifies the connection between $\mathcal{I}(A)$ and
$\mu$ when $A$
is a polygon. Loosely speaking, $\mathcal{I}(A)$ is
simply the $\mu$-length of the polygon.

\begin{lemma}
Let $A$ be a polygon whose sides are $[s_0,s_1],[s_1,s_2],\ldots,[s_{n_e
-1},\break s_{n_e}]$, with $s_{n_e}=s_0$.
We have
\[
\mathcal{I}(A)=\sum_{i=0}^{n_e-1} \mu(s_i-s_{i+1}).
\]
\end{lemma}

\begin{pf}
For each $x=(a,b)\in\Zdeux$, define $x^{\perp}=(-b,a)$.
The map $\R^{\Zdeux}\to\R^{\Zdeux}$ that maps $(t_{x,y})_{\{x,y\}
\in
\Edeux}$
to $(t_{-y,x})_{\{x,y\}\in\Edeux}$ leaves $m^{\otimes\Ed}$
invariant, so
it follows from~(\ref{elementaire}) that $\mu(z)=\mu(z^{\perp})$ holds
for each $z\in\Zdeux$. Since $\mu$ is homogeneous and continuous,
the formula
$\mu(z)=\mu(z^{\perp})$ also holds for each $z\in\Rdeux$.
We have
\begin{eqnarray*}
\int_{\partial^* A}\mu(\nu_A(x))\, d\mathcal{H}^{1}(x) & =&
\sum_{i=0}^{n_e-1} \|s_i-s_{i+1}\|_2\mu\biggl( \biggl(\frac{s_i-s_{i+1}}{\|
s_i-s_{i+1}\|_2}\biggr)^{\perp} \biggr)\\
& =& \sum_{i=0}^{n_e-1} \|s_i-s_{i+1}\|_2\mu\biggl( \frac{s_i-s_{i+1}}{\|
s_i-s_{i+1}\|_2} \biggr)\\
& =& \sum_{i=0}^{n_e-1} \mu(s_i-s_{i+1}).
\end{eqnarray*}
\upqed
\end{pf}

The next property of $\mathcal{I}$ will be decisive in the proof of
lower large deviations. Basically, it says that
the shortest path surrounding a convex polygon is the frontier
of the polygon itself.

\begin{lemma}\label{compare-aire}
Let $A$, $B$ be two polygons with $B\subset A$. We suppose that
$B$ is convex.
Then $\mathcal{I}(B)\le\mathcal{I}(A)$.
\end{lemma}
\begin{pf}
We proceed by induction on the number $n(A,B)$ of vertices of $B$ which
do not belong to $\partial A$. When $n=0$, we just apply the triangle
inequality.
When $n>1$, we build a polygon $A'$ with $B\subset A'\subset A$,
$\mathcal{I}(A')\le\mathcal{I}(A)$ and $n(A',B)<n(A,B)$ as follows:
let $z$ be a vertex of $B$ which is not in $\partial A$. Since $B$ is
convex, there exists an affine map $\phi$ with $\phi(z)=0$ and $\phi
(x)<0$ for $x$ in $B\backslash\{z\}$.
Let $D$ be the connected component of $z$ in $A\cap\{x\in\Rdeux:\phi
(x)\ge0\}$.
$D$~is a polygon which has a side $F$ in $\{x\in\Rdeux:\phi(x)\ge0\}$.
Note $A'=A\backslash D$.
Denote by $s_a$ and $s_b$ the ends of $F$ and define $\mu(F)=\mu(s_b-s_a)$.
We have
$\mathcal{I}(A)=(\mathcal{I}(A')-\mu(F))+(\mathcal{I}(D)-\mu(F))$.
By the triangle inequality $\mu(F)\le\mathcal{I}(D)/2$, so $\mathcal
{I}(A')\le\mathcal{I}(A)$.
\end{pf}

We will also need convenient approximations of a convex set by convex
polygons. This is the goal of the next lemma.

\begin{lemma}
\label{approx-convexe}
Let $A$ be a bounded convex set with $0$ in the interior of $A$.
For each $\epsilon>0$, there exist convex polygons $P$ and $Q$ such that
\[
0\in P\subset A\subset Q\quad\mbox{and}\quad \mathcal{I}(Q)-\epsilon\le
\mathcal{I}(A)\le\mathcal{I}(P)+\epsilon.
\]
\end{lemma}

\begin{pf}
A proof of the existence of $Q$ can be found in Lachand--Robert and
Oudet~\cite{LRO} in a more general setting. The existence of $P$ is simpler:
let $(A_p)_{p\ge1}$ such that:
\begin{itemize}
\item for each $p\ge1$, $A_p$ is a convex polygon,
\item for each $p\ge1$, $A_p\subset A$,
\item$0\in A_p$ for large $p$,
\item$\lim_{p\to+\infty} \lambda^2( A\backslash A_p)=0$.
\end{itemize}
(e.g., take $A_p$ as the convex hull of $x_1,\dots, x_p$, where
$(x_p)_{p\ge1}$ is dense in $\partial A$: this ensures that $\bigcup_{p\ge1} A_p \supset A\backslash\partial A$.)
For fixed $f\in C^1_c(\Rd,\mathcal{W}_{\mu})$, $A\mapsto\int_{A}
\Divc
f \, d\lambda^2(x)$ is continuous with respect to the $L^1$ convergence
of Borel sets, so
$A\mapsto\mathcal{I}(A)$ is lower semicontinuous.
Then $\mathcal{I}(A)\le \liminf_{p\to+\infty} \mathcal
{I}(A_p)$, so
there exists $p\ge1$ with $\mathcal{I}(A)\le\mathcal
{I}(A_p)+\epsilon
$ and $0\in A_p$.
\end{pf}

\section{Upper large deviations}
\label{ladeuxieme}
\begin{thm}\label{lesgrandes}
For each $\epsilon>0,$ there exist constants $C_5,C_6>0$, such that
%
\begin{equation}
\label{une-autre-ineg}
\P\bigl(\mincut(nA,\infty)\ge n\mathcal{I}(A)(1+\epsilon)\bigr)\le C_5\exp
(-C_6 n).
\end{equation}
\end{thm}

The proof naturally falls into three parts:
\begin{enumerate}
\item Approximate $nA$ by a polygon.
\item Parallel outside $nA$ (but close to $nA$) the boundary of the
polygon: it creates a new polygon.
\item Hope that successive vertices of the newly created polygon can be
joined by a path which is short enough and does not enter in $nA$.
\end{enumerate}

Therefore, we need a lemma
that would roughly say that one can find a path from $x$ to $y$
that has length smaller than $(1+\epsilon)\mu(x-y)$ and is not far from
a straight line. To this aim, we introduce some definitions:

Let $y,z\in\Rd$, $\hat x\in\mathcal S$, and $R,h>0$. We define
\begin{eqnarray}
d(y,\R\hat x) & = &\| y-\langle y,\hat x\rangle\hat x\|_2\nonumber\\
\eqntext{\mbox{(the
Euclidean distance from $y$ to the line $\R\hat x$),}} \\
\Cyl_z(\hat x,R,h) & =& \{y\in\Zd\dvtx d(y-z,\R\hat x)\le R\mbox{ and }0\le\langle y-z,\hat x \rangle\le h\},\nonumber
\end{eqnarray}
For $R>0$ and $z,z'\in\Rd$ with $z\ne z'$, we also define
\[
\Cyldeux(z,z',R)=\Cyl_z \biggl(\frac{z'-z}{\|z'-z\|_2},R,\|z-z'\|_2 \biggr).
\]

\begin{lemma}\label{GDlong}
Let $z\in\Rd$, $\hat x\in\mathcal S$, $h\ge1$ and $r\ge1$. We can define
$s_0$ (\resp$s_f$) to be the integer point in $\Cyl_z(\hat x,r,h)$
which is
the closest to $z$ (\resp$z+h\hat x$).
We also define the longitudinal crossing time $\tlong(\Cyl_z(\hat
x,h,r))$ of the cylinder
$\Cyl_z(\hat x,r,h)$ as the minimal time needed to cross it from $s_0$
to $s_f$, using only edges inside the cylinder.

Then for each $\epsilon>0$ and each function $f\dvtx\R_+\to\R_+$ with
$\lim_{r\to+\infty} f(r)=+\infty$,
there exist two strictly positive constants $C_7$ and $C_8$ such that
\begin{eqnarray*}
&&\forall z\in\Rd,\ \forall\hat x \in\mathcal S, \ \forall h>0\\
&&\qquad\P\bigl( \tlong(\Cyl_z(\hat x,f(h),h)) \ge\mu( \hat x) ( 1+ \epsilon) h
\bigr)  \le C_7 \exp(-C_8h).
\end{eqnarray*}
\end{lemma}

\begin{pf}
For $x\in\Zd$ and $t\ge0$, let
\[
\mathcal{B}_{x}(t)=\{y\in\Zd\dvtx \mu(x-y)\le t\}.
\]
For $x,y\in\Zd$ denote by $I_{x,y}$ the length of the shortest path from
$x$ to $y$ which is inside $\mathcal{B}_x(1,25\mu(x-y))\cap\mathcal
{B}_y(1,25\mu(x-y))$. Since $I_{x,y}$ as the same law than $I_{0,x-y}$,
we simply define $I_x=I_{0,x}$.
We begin with an intermediary lemma.
\begin{lemma}
For each $\epsilon\in(0,1]$, there exists $M_0=M_0(\epsilon)$ such that
for each $M\ge M_0$ there exist $c=c(\epsilon,M)<1$ and
$t=t(\epsilon,M)>0$ with
\[
\|x\|\in[M/2,2M] \quad\Longrightarrow\quad\E\exp\bigl(t \bigl(I_x-(1+\epsilon)\mu
(x)\bigr)\bigr)\le c.
\]
\end{lemma}
\begin{pf}
Let $Y$ be a random variable with law $m$ and let $\gamma>0$ be such
that $\E e^{2\gamma Y}<+\infty$.
First, 
equation~(\ref{elementaire})
easily implies the following almost sure convergence:
\[
\lim_{\|x\|\to+\infty} \frac{I_x}{\mu(x)}=1.
\]
By considering a deterministic path from $0$ to $x$ with length $\|x\|
$, we see that $I_x$ is dominated by a sum of $\|x\|$ independent
copies of $Y$ denoted by $Y_1,\dots, Y_{\|x\|}$, and thus
$I_x/{\|x\|}$ is dominated by
\[
\displaystyle\frac1{\|x\|}\sum_{k=1}^{\|x\|} Y_i.
\]
This family is equi-integrable by the law of large numbers. So
$\displaystyle(I_x/{\|x\|})_{x\in\Zd\backslash\{0\}}$ and then
$\displaystyle(I_x/{\mu(x)})_{x\in\Zd\backslash\{0\}}$ are also
equi-integrable families, which implies that
%
\begin{equation}
\label{equiinteg}
\lim_{\|x\|\to+\infty} \frac{\E I_x}{\mu(x)}=1.
\end{equation}
Note that for every $y \in\R$ and $t\in(0,\gamma]$,
\[
e^{ty} \le1+ty +\frac{t^2}2y^2 e^{t|y|}
\le1+ty+\frac{t^2}{\gamma^2} e^{2\gamma|y|}.
\]
Let $\tilde{I}_x=I_x-(1+\epsilon)\mu(x)$ and suppose that $t\in
(0,\gamma]$. Then since $|\tilde{I}_x|\le I_x+2 \mu(x)$, the previous
inequality implies that
\[
e^{t\tilde{I}_x} \le1+t\tilde{I}_x +\frac{t^2}{\gamma^2}
e^{4\gamma\mu
(x)}e^{2\gamma{I}_x}.
\]
Since $\mu(x)\le\|x\|\mumax$ and $I_x \le Y_1+ \cdots+Y_{\|x\|}$,
we can define $\rho=e^{4\gamma\mumax}\E e^{2\gamma Y},$ and thus obtain
\[
\E e^{t\tilde{I}_x}\le1+t\biggl[\E\tilde{I}_x+\frac{t}{\gamma^2}\rho
^{\|x\|}\biggr].
\]
Considering equation~(\ref{equiinteg}), let $M_0$ be such that
$\|x\| \ge M_0/2$ implies $\frac{\E I_x}{\mu(x)}\le1+\epsilon/3$.
For $x$ such that $\|x\|\ge M_0$, we have $\E\tilde{I}_x \le-\frac
23\epsilon\mu(x)$,
so
\[
\E e^{t\tilde{I}_x} \le1+t\biggl[ -\frac23\epsilon\mu(x)+\frac
{t}{\gamma
^2}\rho^{\|x\|}\biggr]
\le1+t\biggl[ -\frac23\epsilon\mumax+\frac{t}{\gamma^2}\rho^{\|x\|}\biggr].
\]
Therefore, we can take $t=t(\epsilon,M)=\min(\gamma,\gamma^2\mumax
\frac
{\epsilon}3\rho^{-2M})$ and $c=c(\epsilon,M)=1-\frac13\epsilon
\mumax
t(\epsilon,M)$.
\end{pf}

Let us come back now to the proof of Lemma~\ref{GDlong}.
Let $\epsilon\in(0,1)$ and consider the integer $M_0=M_0(\epsilon/3)$
given by the previous lemma.
Let $M_1=M_1(\epsilon)$ be an integer greater than $M_0$ and such that
%
\begin{equation}
\label{condM}
(1+\epsilon/3) \biggl(1+\frac{\mumax}{\mumin}\frac{2}{M_1} \biggr)\le
1+\epsilon/2.
\end{equation}
%
Let $N$ be the smallest integer which is greater than $h/M_1$ and, for
each $i\in\{0,\dots,N\},$
denote by $x_i$ the integer point in the cylinder which is the closest to
$z+\frac{ih\hat x}{N}$. Note that
\[
\biggl(1-\frac1N \biggr) M_1 \le\frac h N \le M_1.
\]

1. Let $i_0$ be an integer with $i_0\ge\max(\frac{1,25(2+M_1)\sqrt
2}{\mumin},2 )$.\vspace*{2pt}
There exists a deterministic path inside the cylinder from $x_0$ to $x_{i_0}$
(\resp$x_{N-i_0}$ to $x_N$) which uses less than $2i_0h/N$ edges: we
denote by
$L_{\mathrm{start}}$ (\resp$L_{\mathrm{end}}$) the random length of this path.
Markov's inequality easily gives
%
\begin{eqnarray}
\label{lepetitbout}
& &\P\biggl(L_{\mathrm{start}}>\frac{\epsilon}4 \mu(\hat x)h\biggr)+\P\biggl(L_{\mathrm{end}}>\frac
{\epsilon
}4 \mu(\hat x)h\biggr)\\
&&\qquad\le 2(\E e^{2\gamma Y})^{2i_0M_1}\exp\biggl(-\frac
{\gamma\epsilon}2 h\mumin\biggr)\le C'e^{-C_5 h}.
\end{eqnarray}

2. For each $i, j\in\{0,\dots,N-1\}$, we have $|\mu
(x_i-x_{j})-|j-i|h\mu
(\hat x)/N|\le2\mumax$. Thus, if $h$ is larger than some $h_0$, then
$\mathcal{B}_{x_i}(1,25\mu(x_i-x_{i+1})) \cap\mathcal
{B}_{x_j}(1,25\mu
(x_j-x_{j+1}))=\varnothing$ as soon as $|j-i|\ge2$.

Let $h_1=h_1(\epsilon,f) \ge h_0$ be such that
$\forall h \ge h_0, f(h)\ge i_0$. 
If we take $h$ larger than $h_1$, then the whole set
\[
\bigcup^{N-i_0-1}_{i=i_0}\mathcal{B}_{x_i}\bigl(1,25\mu(x_i-x_{i+1})\bigr)
\]
stays
inside the cylinder.
So, provided that $h\ge h_1$, we have inside the cylinder a path from
$x_0$ to $x_N$
with length
\[
L_{\mathrm{start}}+\sum_{i=i_0}^{N-i_0-1}
I_{x_i,x_{i+1}}+L_{\mathrm{end}}.
\]
%
Let
\[
S_{\mathrm{odd}}=\mathop{\sum_{2\le i\le N-3}}_{ i\ \mathrm{odd}} I_{x_i,x_{i+1}}
\quad\mbox{and}\quad
S_{\mathrm{even}}=\mathop{\sum_{2\le i\le N-3}}_{ i\ \mathrm{even}} I_{x_i,x_{i+1}}.
\]

By the definition of $(x_i)_{1\le i\le N}$, we have
\begin{eqnarray*}
\mathop{\sum_{2\le i\le N-3}}_{ i\ \mathrm{odd}} \mu
(x_{i+1}-x_i) & \le&\mathop{\sum_{2\le i\le N-3}}_{ i\ \mathrm{odd}}
\frac{h\mu(\hat x)}{N}+2\mumax\\
& \le&\frac{N-3}2 \biggl(\frac{h\mu(\hat x)}{N}+2\mumax\biggr)\\
& \le&\frac{N}2\frac{h\mu(\hat x)}{N}+ (N-1)\frac\mumax\mumin\mu
(\hat
x)\\
& \le&\frac{h\mu(\hat x)}2\biggl(1+2\frac\mumax\mumin\frac1{M_1}\biggr)
\end{eqnarray*}

Then using~(\ref{condM}), we can write, for each $t\ge0$,
\begin{eqnarray*}
&&\P\biggl(S_{\mathrm{odd}}\ge\frac{h\mu(\hat x)}2(1+\epsilon/2)\biggr) \\
&&\qquad \le\P
\Biggl(S_{\mathrm{odd}}\ge
(1+\epsilon/3) \mathop{\sum_{2\le i\le N-3}}_{ i\ \mathrm{odd}} \mu
(x_{i+1}-x_i)\Biggr)\\
&&\qquad \le\P\Biggl( \mathop{\sum_{2\le i\le N-3}}_{ i\ \mathrm{odd}}
I_{x_i,x_{i+1}}-(1+\epsilon/3)\mu(x_{i+1}-x_i)\ge0\Biggr)\\
&&\qquad \le\P\Biggl( \exp\Biggl(t\mathop{\sum_{2\le i\le N-3}}_{ i\ \mathrm{odd}}
I_{x_i,x_{i+1}}-(1+\epsilon/3)\mu(x_{i+1}-x_i)\Biggr)\ge1\Biggr)\\
&&\qquad \le\E\exp\Biggl(t\mathop{\sum_{2\le i\le N-3}}_{ i\ \mathrm{odd}}
I_{x_i,x_{i+1}}-(1+\epsilon/3)\mu(x_{i+1}-x_i)\Biggr)\\
&&\qquad \le\mathop{\prod_{2\le i\le N-3}}_{ i\ \mathrm{odd}} \E\exp
\bigl(tI_{x_i,x_{i+1}}-(1+\epsilon/3)\mu(x_{i+1}-x_i)\bigr).
\end{eqnarray*}
We take now $t=t(\epsilon/3,M_1(\epsilon))$ and $\rho=\rho(\epsilon
/3,M_1(\epsilon))$. For each $i$, we have $\mu(x_i-x_{i+1})\in
[M_1/2,2M_1]$, thus we can apply the previous lemma and get
\begin{eqnarray*}
\P\biggl(S_{\mathrm{odd}}\ge\frac{h\mu(\hat x)}2(1+\epsilon/2)\biggr) & \le&\rho
^{(N-5)/2}\\
& \le&\rho^{hM_1(\epsilon)/2-3/2}=A\exp(-Bh)
\end{eqnarray*}
with $A=\rho^{-3/2}$ and $B=-\frac1{2M_1(\epsilon)}\ln\rho$.

Similarly, $\P(S_{\mathrm{even}}\ge\frac{h\mu(\hat x)}2(1+\epsilon/2))\le
A\exp
(-Bh)$, so it suffices to put the pieces together to conclude the proof.
\end{pf}

\begin{pf*}{Proof of Theorem~\ref{lesgrandes}}
We first consider the case where $A$ is a convex polygon.
Let us denote by $s_0,s_1,\dots,s_{n_e}$ the vertices of $A$, with
$s_{n_e}=s_0$.
We suppose that the vertices are in trigonometric order.
For each $i\in\{0,n_e-1\}$, let $v_i$ be such that
$\langle v_i,s_{i+1}-s_i\rangle=0$ and $\langle v_i,s_i\rangle=1$.
For $x\in\Rd$, define $\phi_i(x)=\langle v_i,x\rangle$. With our conventions
\[
nA=\bigcap^{n_e-1}_{i=0}\{x\in\Rdeux\dvtx \phi_i(x)\le n\}.
\]

For $z\in\Rdeux$, we define $\Int(z)$ as the only $x\in\Zp$ such that
$z\in x+[-1/2,1/2)\times[-1/2,1/2)$. Let $\epsilon>0$.
For $i\in\{0,\ldots, n_e\}$, let $y_i=\Int(n(1+\epsilon)s_i )$.

Our goal is to build for each $i$ a path from $y_i$ to $y_{i+1}$
which does not enter $nA$ and is short enough.
Define $M=\max\{\|v_i\|_2;0\le i\le n_e-1\}$
and $S=\max\{\mu(s_i-s_{i+1});0\le i\le n_e-1\}$.

It is easy to see that
\[
\forall r\ge0\qquad\phi_i\ge n(1+\epsilon)-M\biggl(\frac{\sqrt
2}2+r\biggr)\qquad{on }\
\Cyldeux(y_i,y_{i+1},r).
\]

Moreover, for each $i\in\{0,\dots,n_e-1\}$, we have
\[
M\biggl(\frac{\sqrt2}2+\frac{\epsilon}{4\mathit{MS}}\|y_i-y_{i+1}\|_2\biggr)
\le M\biggl(\frac{\sqrt2}2+\frac{\epsilon}{4\mathit{MS}}\bigl(nS+\sqrt2\bigr)\biggr)\le\frac
{n\epsilon}2,
\]
provided that $n$ is large enough. Therefore, it follows
that $\phi_i\ge(1+\epsilon/2)n$ on 
$\displaystyle\Cyldeux(y_i,y_{i+1},\frac{\epsilon}{4\mathit{MS}}\|
y_i-y_{i+1}\|
_2)$, which means that this set is off $nA$.

Since $\mu(y_i-(1+\epsilon)n s_i)\le\mumax$, we know that
\[
\sum_{i\in\{0,\dots, n_e\}}\mu(y_i-y_{i+1})\le
n(1+\epsilon
)\mathcal{I}(A)+2n_e\mumax\le n(1+\epsilon)^2\mathcal{I}(A)
\]
provided that $n$ is large enough.

Then one can see that for $n$ greater than some (deterministic) integer
$n_0$, the event
\begin{eqnarray*}
A_n &=& \bigcap_{i\in\{0,\dots, n_e-1\}} \biggl\{
\tlong\biggl(
\Cyldeux\biggl(y_i,y_{i+1},\frac{\epsilon}{4\mathit{MS}}\|y_i-y_{i+1}\|_2\biggr)
\biggr)\\
&&\hspace*{140pt}
<(1+\epsilon)\mu(y_{i+1}-y_i) \biggr\}
\end{eqnarray*}
satisfies
\[
A_n \subset\{\Cut(nA)\le n(1+\epsilon)^3\mathcal{I}(A)\}
\]

We are now ready to apply Lemma~\ref{GDlong} with
$f(h)=\frac{\epsilon}{4\mathit{MS}}h$.
It comes that
\begin{eqnarray*}
P\bigl(\Cut(nA)> n(1+\epsilon)^3\mathcal{I}(A)\bigr) & \le&
\P(A_n^c)\\
& \le& \sum^{n_e-1}_{i=0} C_7\exp(-C_8\|y_i-y_{i+1}\|_2)\\
& \le&\sum^{n_e-1}_{i=0} C_7e^{C_8\sqrt2}\exp(-C_8 \|
s_i-s_{i+1}\|_2 n),\\
& \le& c_1 \exp(-c_2 n),
\end{eqnarray*}
%
with $c_1= n_eC_7e^{C_8\sqrt2}$ and $c_2=C_8\min_i \|s_i-s_{i+1}\|_2$.

Since $\epsilon$ is arbitrary, the theorem follows when $A$ is a polygon.

Let us go to the general case: By Lemma~\ref{approx-convexe}, there
exists a convex polygon
$Q$ with $Q\supset A$ and
$(1+\epsilon)\mathcal{I}(A)\le(1+\epsilon/2) \mathcal{I}(Q).$

By its very definition, $\mincut(nA,\infty)\le\mincut(nQ,\infty)$.
Then
\begin{eqnarray*}
\P\bigl(\mincut(nA,\infty)\ge n\mathcal{I}(A)(1+\epsilon)\bigr) & \le&\P
\bigl(\mincut
(nQ,\infty)\ge n\mathcal{I}(A)(1+\epsilon)\bigr)\\
& \le&\P\bigl(\mincut(nQ,\infty)\le n\mathcal{I}(Q)(1+\epsilon/2)\bigr).
\end{eqnarray*}

Hence, the result follows from the polygonal case.
\end{pf*}

\section{Lower large deviations}
\label{latroisieme}

\begin{thm}\label{lespetites}
For each $\epsilon>0$,
there exist constants $C_9,C_{10}>0$, such that
%
\begin{equation}\label{une-autre-ineg-dessous}
\forall n\ge1\qquad\P\bigl(\mincut(nA,\infty)\le n\mathcal
{I}(A)(1-\epsilon)\bigr)\le C_9\exp(-C_{10} n).
\end{equation}
\end{thm}

The choice of a strategy for the proof of lower large deviations is more
difficult than for the upper ones. An important point is that it is
hopeless to consider the sides of the polygon separately.

\begin{figure}

\includegraphics{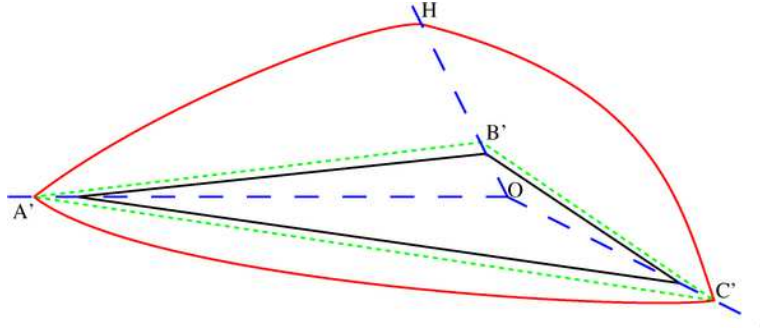}

\caption{Surrounding the polygon.}
\label{unefigure}
\end{figure}

Indeed, consider the following picture on Figure~\ref{unefigure}: the
red curve and the green one surround
the black triangle. Of course, it is expected that the minimal cutset looks
like the green triangle rather than like the red ones.
However, the red path from $A'$ to $H$ is shorter than the green one
from $A'$
to $B'$. But this advantage is lost on the next side because
the red path from $H$ to $C'$ is much longer than the green one from $B'$
to $C'$. So, it appears that we must think globally, using the perimeter
of surrounding curves. To this aim, Lemma~\ref{compare-aire} will be
particularly useful.

\begin{pf*}{Proof of Theorem \ref{lespetites}}
Again, we first deal with the case, where $A$ is a convex polygon
whose sides are $[s_0,s_1],[s_1,s_2],\dots,[s_{n_e-1},s_{n_e}]$, with
$s_{n_e}=s_0$.
We denote by $L_{n,i}$ the points $x\in\Zp$ that touch a bond which intersects
$[1,+\infty)n s_i$.

\begin{lemma}
\label{touslespoints}
For each $\epsilon>0$, there exist $C_{11}=C_{11}(\epsilon
),C_{12}=C_{12}(\epsilon)$, such that
\begin{eqnarray*}
&&\P\bigl(\exists i\in\{0,n_e-1\} \exists(x,y)\in L_{n,i}\times
L_{n,i+1}\, d(x,y) \le (1-\epsilon)\mu(x-y) \bigr)\\
&&\qquad\le C_{11}\exp
(-C_{12} n).
\end{eqnarray*}
\end{lemma}

\begin{pf}
Since $\{0,\dots,n_e-1\}$ is finite, it is sufficient to prove that
for each
$i,j$ with $0\le i< j<n_e$, there exists $C_{11}(i,j)>0$ and
$C_{12}(i,j)>0$ with
\[
\P\bigl( \exists(x,y)\in L_{n,i}\times L_{n,j} \,d(x,y)\le(1-\epsilon
)\mu(x-y) \bigr)\le C_{11}(i,j)\exp(-C_{12}(i,j) n).
\]

Thanks to Proposition~\ref{mugs}, we can write
\begin{eqnarray*}
&&\P\bigl(\exists(x,y)\in L_{n,i}\times L_{n,j} d(x,y)\le
(1-\epsilon)\mu(x-y) \bigr)\\
&&\qquad\le \sum_{(x,y)\in L_{n,i}\times L_{n,j}} \P\bigl(d(x,y)\le
(1-\epsilon)\mu(x-y) \bigr)\\
&&\qquad \le \sum_{(x,y)\in L_{n,i}\times L_{n,j}} C_3\exp
(-C_4\|x-y\|_2)\\
&&\qquad \le C_3 \sum^{+\infty}_{p=0} |A_p| \exp(-C_4 p),
\end{eqnarray*}
where
\[
A_p=\{(x,y)\in L_{n,i}\times L_{n,j}; \|x-y\|_2\in[p,p+1)\}.
\]
Let $\alpha=d_2([1,+\infty)s_i,[1,+\infty)s_j)$
and $\theta=\arccos\frac{\langle s_i,s_j\rangle}{\|s_i\|_2\|s_j\|_2}$.
We can see that:
\begin{itemize}
\item$|A_p|=0$ for $p\le n\alpha-3$.
\item$|A_p|\le\frac{2000}{\sin^2\theta} (1+p)^2$ for each $p\ge0$.
\end{itemize}
The first point is clear. Let us prove the second point: for each $k\in
\{i,j\}$, let $s'_k=s_k/\|s_k\|_2$.
Obviously, $A_p\subset B_p\times B'_p$, where
$B_p=\{x\in L_{n,i} ; d_2(x,\R s'_j)\le p+3\}$ and
$B'_p=\{y\in L_{n,j} ; d_2(x,\R s'_i)\le p+3\}$.

For $r\in\R$, define
\[
f(r)=\sum_{x\in B_p} \1_{\{|r s'_i-x|\le\sqrt2\}}.
\]
Since $d_2(x,\R s_i)\le1$ for each $x\in L_{n,i}$, it follows that
\[
\displaystyle\int_{\R} f(r)\,dr\ge2|B_p|.
\]
For a given $r$, the sum defining $f(r)$ has at most 9 nonvanishing
terms, thus we have
\[
f(r)\le9 \1_{\{d_2(r s'_i,\R s_j)\le p+3+\sqrt2\}}.
\]
Then
\[
|B_p|\le\frac12\int_{\R} f(r)\,dr \le9\times\frac{1}{\beta
_{i,j}}\bigl(p+3+\sqrt2\bigr)\le\frac{9(3+\sqrt2)}{\beta_{i,j}}(p+1),
\]
where $\beta_{i,j}=|s'_i-\langle s'_i,s'_j\rangle s'_j|=\sqrt
{1-\langle
s'_i,s'_j\rangle^2}$.\\
Similarly, $|B'_p|\le\frac{9(3+\sqrt2)}{\beta_{i,j}}(p+1)$.
Finally, $|A_p|\le\frac{2000}{\sin^2 \theta}(1+p)^2$.

Let $K'$ be such that $ \frac{2000}{\sin^2 \theta} C_3(1+p)^2\le
K'\exp
(\frac{C_4}{2} p)$ holds for each $p\ge0$:
we have
\begin{eqnarray*}
& & \P\bigl(\exists(x,y)\in L_{n,i}\times L_{n,j}\dvtx d(x,y)\le
(1-\epsilon)\mu(x-y)\bigr)\\
&&\qquad \le \sum^{+\infty}_{p=\Int( n\alpha-3)} K'\exp\biggl(-\frac
{C_4}{2} p\biggr)
\le\frac{K'e^{2 C_4}}{1-\exp(-C_4/2)}\exp\biggl(-\frac
{C_4\alpha}2 n\biggr),
\end{eqnarray*}
which completes the proof of the lemma.
\end{pf}

We go back to the proof of Theorem~\ref{lespetites}.

Suppose that $\mincut(nA,\infty)<(1-\epsilon)n\mathcal{I}(A)$. Then we
can find in the dual lattice a closed path $\gamma$ that surrounds $nA$
and whose length $l(\gamma)$ is smaller than $(1-\epsilon)n\mathcal
{I}(A)$. $\gamma$ necessarily cuts the half-lines $([1,+\infty
)ns_i)_{0\le i\le n_e-1}$ in some points $y_0,y_1,\dots, y_{n_e-1}$. We
also define $y_e=y_0$. The points can be numbered
in such a way that $\gamma$ visits the $(y_i)_{0\le i\le n_e}$ in the
natural order.
Let $x_i$ the point in $L_{n,i}$ which is such that $\|y_i-x_i\|_1\le1/2$.
Obviously,
%
\begin{equation}
\label{petited}
\sum_{i=0}^{n_e-1} d(x_i,x_{i+1})\le l(\gamma)\le(1-\epsilon
)n\mathcal{I}(A).
\end{equation}

Let $B$ be the polygon determined by the $y_i$: we
have
\[
\mathcal{I}(B) =  \sum_{i=0}^{n_e-1}\mu(y_i-y_{i+1})
\le
\sum_{i=0}^{n_e-1} \bigl(\mu(x_i-x_{i+1}) +\mumax\bigr)
\]
$nA$ is convex and contained in $B$, so by Lemma~\ref{compare-aire},
$\mathcal{I}(B)\ge\mathcal{I}(nA)$. It follows that
\begin{eqnarray*}
\sum_{i=0}^{n_e-1} d(x_i,x_{i+1})& \le&(1-\epsilon)n\mathcal{I}(A)
\le
(1-\epsilon) \mathcal{I}(B)\\
& \le&\sum_{i=0}^{n_e-1} (1-\epsilon) \bigl(\mu(x_i-x_{i+1}) +\mumax\bigr)\\
& \le&\sum_{i=0}^{n_e-1} (1-\epsilon/2)\mu(x_i-x_{i+1}),
\end{eqnarray*}
provided that $n\ge\frac1\alpha(1+\frac2\epsilon\frac\mumax
\mumin)$.

So, for large $n$, the event $\{\mincut(nA,\infty)<(1-\epsilon)\}$
implies the existence of $i\in\{0,\dots,n_e-1\}$, $x_i\in L_{n,i}$ and
$x_{i+1}\in L_{n,i+1}$ with
\[
d(x_i,x_{i+1})\le(1-\epsilon/2)\mu(x_i-x_{i+1}).
\]
Then we have
\begin{eqnarray*}
& & \P\bigl(\mincut(nA,\infty)<(1-\epsilon)n\mathcal
{I}(A)\bigr) \\
&&\qquad \le \P\bigl(\exists i\in\{0,\dots,\ne-1\},\\
&&\qquad\phantom{\le \P\bigl(}
{}\exists(x,y)\in
L_{n,i}\times L_{n,i+1}\, d(x,y)\le(1-\epsilon/2)\mu(x-y)\bigr)\\
&&\qquad \le C_{11}(\epsilon/2)\exp(-C_{12}(\epsilon/2) n),
\end{eqnarray*}
thanks to Lemma~\ref{touslespoints}. This ends the proof in the case,
where $A$ is a polygon.

Let us go to the general case: By Lemma~\ref{approx-convexe}, there
exists a convex polygon
$P$ with $0\in P$, $P\subset A$ and
$(1-\epsilon)\mathcal{I}(A)\le(1-\epsilon/2) \mathcal{I}(P).$

By its very definition, $\mincut(nA,\infty)\ge\mincut(nP)$.
Then
\begin{eqnarray*}
&&\P\bigl(\mincut(nA,\infty)\le n\mathcal{I}(A)(1-\epsilon)\bigr) \\
&&\qquad \le\P
\bigl(\mincut
(nP)\le n\mathcal{I}(A)(1-\epsilon)\bigr)\\
&&\qquad \le\P\bigl(\mincut(nP)\le n\mathcal{I}(P)(1-\epsilon/2)\bigr),
\end{eqnarray*}
which has just been proved to decrease exponentially fast with
$n$.\
\end{pf*}

\section{Final proofs}
\label{laquatrieme}
\subsection{Proof of the theorems}

Obviously, Theorems~\ref{lesgrandes} and \ref{lespetites} concur
to get Theorem~\ref{vitesseexpo}.
Since $\mathcal{I}(A)=\int_{\partial^* A}\mu(\nu_A(x))\, d\mathcal
{H}^{1}(x)$,
the first equality in Theorem~\ref{main} directly follows from
Theorem~\ref{vitesseexpo} with the help of the Borel--Cantelli lemma.

It is worth saying a word about equation~(\ref{maxflowmincutinfini}), because
the Ford--Fulkerson theorem is initially concerned with finite graphs.
Let us recall a version of this theorem.\vspace*{-2pt}
\begin{prop}[(Ford--Fulkerson)]
For each finite graph $G=(V,E)$ and every disjoint subsets $A$ and $B$
of $V$, we have
%
\begin{equation}
\label{maxflowmincutfini}
\max\Biggl\{ \sum_{x\in A}\Div j(x); j\in\Capflow(A,B)\Biggr\}=\mincut(A,B),\vspace*{-2pt}
\end{equation}
where
%
\begin{equation}
\quad \mincut(A,B)=\min\Biggl\{\sum_{x\in C} t_x; \mbox{ every path in $G$ from $A$
to $B$ meets $C$}\Biggr\}\hspace*{-12pt}
\end{equation}
and $\Capflow(A,B)$ is the set of flows $j$ that satisfy $|j(x,y)|\le
t_{\{x,y\}}$ for each $\{x,y\}\in E$ and $\Div j (x)=0$\ for $x\in
V\backslash(A\cup B)$.
\end{prop}

In fact, in the initial paper~\cite{MR0079251} and in most books, $A$
and $B$ are just singletons. The reduction to this case is easy.
Because of the antisymmetry property, the contribution of edges inside
$A$ to $\sum_{x\in A}\Div j(x)$ is null; so we neither change the
max-flow nor the min-cut if we identify
the points that are in $A$. Obviously, the max-flow and the min-cut are
not changed either when we identify the points that are in
$B$.

Now let $G_n=(V_n,E_n)$ be the restriction of $\Ld$ to $V_n=\{x\in\Zd;
\|x\|_1\le n\}$ and denote by $B_n$ the boundary of $V_n$.

Let $f$ be a flow from $A$ to infinity; particularly, $f$ is a flow
from $A$ to $B_n$,
so $\sum_{x\in A}\Div j(x)\le\mincut(A,B_n)$.
By the definition of a cutset, a minimal cutset from $A$ to infinity is the
external boundary of a finite connected set containing~$A$.
In particular, a minimal cutset is finite.
It follows that $\inf_{n\ge1} \mincut(A,B_n)=\mincut(A,\infty)$.
Then $\sup\{ \sum_{x\in A}\Div j(x); \break j\in\Capflow(A,\infty)\}\le
\mincut
(A,\infty)$.
Conversely, let $j_n$ be a flow that realizes
$\max\{ \sum_{x\in A}\Div j(x); j\in\Capflow(A,B_n)\}$.
We can extend $j_n$ to $\Edo$ by putting $j_n(e)=0$ outside $E_n$.
Obviously, $j_n\in\prod_{e\in\Edo} [-t_e,+t_e]$, thus the sequence
$(j_n)_{n\ge1}$ admits a subsequence $(j_{n_k})_{k\ge1}$ converging to
some $j'\in\prod_{e\in\Edo} [-t_e,+t_e]$ in the product topology.
Easily, $j'$ is antisymmetric.\vspace*{-2pt}
%
\begin{eqnarray*}
\sum_{x\in A}\Div j'(x) & = &\lim_{k\to+\infty} \sum_{x\in A}\Div
j_{n_k}(x)\\
& =& \lim_{k\to+\infty} \max\Biggl\{ \sum_{x\in A}\Div j(x); j\in
\Capflow
(A,B_k)\Biggr\}\\
& =& \lim_{k\to+\infty} \mincut(A,B_{n_k})\\
& =& \inf_{n\ge1} \mincut(A,B_{n})= \mincut(A,\infty).
\end{eqnarray*}
For each $x\in\Zd\backslash A$, there exists $k_0$ such that $x\in
V_n\backslash( B_{n_k}\cup A)$ for $k\ge k_0$; then $\Div j_{n_k}(x)=0$
for $k\ge k_0$, which ensures that $\Div j'(x)=0$. It is now easy to
see that $j'$ is a capacitive flow from $A$ to infinity, which
completes the proof of equation~(\ref{maxflowmincutinfini}) and,
therefore, the proof of Theorem~\ref{main}.

\subsection{Proof of the corollaries}

Let us now recall Menger's theorem (see, for instance, Diestel~\cite{MR2159259} for a proof).
\begin{prop}[(Menger's theorem)]
Let $G=(V,E)$ be a finite graph and $A,B\subset V$. Then the minimum
number of vertices
separating $A$ from $B$ is equal to the maximum number of disjoint
paths from $A$ to $B$.
\end{prop}

We can now prove Corollary~\ref{avecmenger}.

\begin{pf*}{Proof of Corollary \ref{avecmenger}}
Consider the probability space $(\Omega,\mathcal{B},\P)$, with
$\Omega=\{0,1\}^{\Ed}$ and $\P=\Ber(p)^{\otimes\Ed}$.
As usual, $e$ is said to be open if $\omega_e=1$ and closed otherwise.
Let $R=\{e\in\Ed\dvtx \omega_e=1\}$ and define $V_n$ and $E_n$ as previously.
Let $H_n=(V_n,E_n\cap R)$. It is easy to see that the minimum number of vertices
separating $A$ from $B_n$ is equal to $\mincut(A,B_{n})$, where the
capacity flow is defined by $t_e=1-\omega_e$.
Then by Menger's theorem, the maximum number of disjoint paths from $A$ to
$B_n$ is $\mincut(A,B_{n})$.
By a classical compactness argument, the maximum number of disjoint
paths from $A$ to infinity is the limit of the maximum number of
disjoint paths from $A$ to
$B_n$. Therefore, $\dis(A)=\lim_{n\to+\infty} \mincut
(A,B_{n})=\mincut(A,\infty)$.
The variables $(t_e)_{e\in\Ed}$ are independent Bernoulli variables
with parameter $1-p$. Note $m$ for their common distribution.
Since $p>p_c(2)=1/2$, $m(0)=1-p<1/2$, and we can apply Theorem~\ref
{vitesseexpo} to complete the proof of Corollary~\ref{avecmenger}.
\end{pf*}

We finally prove Corollary~\ref{avecmenger2}.
\begin{pf*}{Proof of Corollary~\ref{avecmenger2}}
Let us denote by $I_k$ the event: ``there exist $k$ disjoint open
biinfinite paths.'' $I_k$ is obviously translation-invariant, so by the
ergodic theorem,
its probability is null or full.
Let $A=[-1,1]^2$ and $S_n=\{\dis(An)\ge n\mathcal{I}(A)/2\}$.
For large $n$, we have $n\mathcal{I}(A)/2>2k$ and
$\P(S_n)>1/2$.
Now consider the event $T_n$: ``all edges inside $nA$ are open.''
It is not difficult to see that $T_n\cap S_n\subset I_k$
but $T_n$ and $S_n$ are independent, so $P(I_k)\ge P(T_n\cap
S_n)=P(T_n)P(S_n)>0$. Finally, $P(I_k)=1$.
\end{pf*}

\section{Perspectives}

It is to be expected that these results still hold in higher
dimensions. In fact, we make the following conjecture:

\begin{conjecture}
We suppose that $m(0)<1-p_c(\Z^d)$ and that
\[
\int_{[0, \infty)}\exp(cx) \,d m(x)<+\infty
\]
for some $c>0$.
Then there exists a map $\mu$ on the unit sphere such that
for each convex set $A$ with 0 in the interior, we have
\[
\lim_{n\to+\infty} \frac{\mincut(nA,\infty)}{n^{d-1}}=\int
_{\partial^*
A}\mu(\nu_A(x)) \,d\mathcal{H}^{d-1}(x).
\]
%
\end{conjecture}

Equivalently,
\begin{eqnarray*}
& &\lim_{n\to+\infty} \frac1{n^{d-1}}\max\Biggl\{ \sum_{x\in nA\cap\Z
^d}\Div j(x); j\in\Capflow(nA,\infty) \Biggr\}\\
&&\qquad= \sup\biggl\{\int_{A}
\Divc f
 \,d\lambda^d(x);f\in C^1_c(\R^d,\mathcal{W}_{\mu}) \biggr\},
\end{eqnarray*}
where
\[
\mathcal{W}_{\mu}=\{x\in\R^d\dvtx \langle x,w\rangle\le\mu(w)\mbox{ for
all }w\}.
\]
Of course, the situation is more complicated when $d\ge3$
because cutsets are not paths; therefore, the capacities can not be interpreted
in term of first-passage percolation.
In a seminal paper~\cite{kesten-surfaces}, Kesten put the basis of a
generalization of
first-passage percolation which seems to be the appropriate tool for the
problem considered here. Basically, he studies the minimal cut between
opposite sides of a parallelepiped with $(e_1,e_2,e_3)$ as axes.
This allows to define a quantity $\nu$ which is a good candidate for
$\mu(e_1)$. Later, Boivin~\cite{boivin-surfaces} extended some of
Kesten's results.
Particularly, he defined a function on the unit sphere of $\R^3$ which
may be
convenient for our purpose.
The condition $m(0)<1-p_c$ is coherent with some previous results; indeed,
Zhang~\cite{MR1749233} proved that $\nu=0$ for $m(0)\ge1-p_c$ whereas
Chayes and Chayes~\cite{MR847132} had proved (at least in the
Bernoulli case) that $\nu>0$ for $m(0)\ge1-p_c$ using a
result of Aizenman, Chayes, Chayes, Fr{\"o}hlich and Russo~\cite{MR728447}. Note that Th\'{e}ret~\cite{theret-small} recently proved
some results that give an independent proof of this fact.
So, $m(0)<1-p_c$ seems to be a natural assumption
for the conjecture. This is also coherent with the expected domain of validity
for the $d$-dimensional version of Corollary~\ref{avecmenger}.
Of course, this conjecture is at present far from being solved because some
of the quantities that are used in the present proof do not have
an obvious equivalent in higher dimensions.
However, we think that the conjecture presented here is a good motivation
to continue the study initiated in Kesten~\cite{kesten-surfaces}.

%

\printaddresses

\end{document}